\documentclass[11pt]{article}

\usepackage{amssymb}

\newcommand{\R}{{\mathbb R}}

\newcommand{\beqn}{\begin{eqnarray*}}
\newcommand{\eeqn}{\end{eqnarray*}}
\newcommand{\beq}{\begin{eqnarray}}
\newcommand{\eeq}{\end{eqnarray}}
\newcommand{\bi}{\begin{itemize}}
\newcommand{\ei}{\end{itemize}}
\newcommand{\text}[1]{\hbox{\rm \ #1\ \/}}
\newcommand{\be}[1]{\begin{equation}\label{#1}}
\newcommand{\ee}{\end{equation}}
\newcommand{\barx}{\overline{x}}
\newcommand{\barz}{\overline{z}}

\title{A Remark on Monotone I/O Systems} 
\author{Eduardo D. Sontag\footnote{Email: sontag@math.rutgers.edu}\\
Dept.\  of Mathematics, Rutgers University, New Brunswick, NJ}

\begin{document}

\maketitle

\begin{abstract}
\noindent
This note remarks that small-gain results for a negative feedback loop around 
a monotone system can be seen as consequences of results concerning an extended
monotone system.
\end{abstract}

\section*{Introduction}

Several recent papers, see e.g.
\cite{angeli-TAC,angeli-SCL,german-DCDS,german-SCL,patrick-MBE,lsu-book}
and references there, have dealt with input/output monotone systems.
Restricting for concreteness to the finite-dimensional continuous-time case,
one considers systems with inputs and outputs in the standard sense of
control theory (\cite{mct}):
\be{sys1}
\dot x = f(x,u),\quad y=h(x)
\ee
and the interest is to study the effect of unity feedback $u=y$,
assuming that the input and output value spaces are the same,
(More general feedbacks can be considered as well, simply redefining $h$.
Dynamic feedback can be included, too, by enlarging the system.)

The general assumption is that $\dot x=f(x,u)$ is monotone as a system with
inputs, meaning that partial orders are given on states and input values,
and the dynamics preserves the orders.
It is also assumed that the open loop system is stable in an appropriate sense.
The state space is a suitable subset of $\R^n$.

There are basically two very distinct types of theorems:
\bi
\item
{\em Positive-feedback results, for which $h$ is a monotone map.}
In this case, one relates location and stability of steady states to
a reduced-dimension discrete or continuous-time dynamics, a dynamics which
evolves in the space of input (and output) values.
\item
{\em Negative-feedback results, for which $h$ is anti-monotone.}
Here one gives sufficient conditions for stability.
The conditions can be interpreted as saying that ``small enough gain'' in the
feedback loop does not destroy open-loop stability. 
\ei

The purpose of this communication is to point out that these two
types of theorems are very closely related.  In summary, studying
the system~(\ref{sys1}) under a negative feedback can be reduced, in the
context of small-gain results, to the study of the following $2n$-dimensional
extended system:
\beq
\dot x&=&f(x,u)\nonumber\\
\label{sys2}
\dot z&=&f(z,h(x))\\
y&=&h(z)\nonumber
\eeq
(obtained by cascading two copies of the original system)
under feedback $u=y$.  The key observation is that 
{\em the extended system is monotone in the i/o sense,\/} 
and therefore the positive-feedback results can be applied to
the closed-loop under $u=y$.  When interpreted in terms of the ``diagonal''
system obtained by restricting to $x=z$, one recovers precisely the small-gain
results.

We will keep the discussion somewhat informal, in order to emphasize the main
ideas.  Missing details can be filled-in by using the techniques developed in
the references.

\section*{Definitions and Assumptions}

We assume existence and uniqueness of maximal solutions of $\dot x=f(x,u)$
under locally bounded inputs; local Lipschitz conditions as in~\cite{mct}
suffice for this.  We take $h$ to be continuous.  
States $x(t)$ are assumed to belong to a closed invariant
convex subset $X$ of $\R^n$, and inputs to a closed convex subset $U$ of
$\R^m$, for some $n,m$.  Thus, solutions starting at $x(0)\in X$ and with inputs
taking values $u(t)\in U$ have the property that $x(t)\in X$ for all $t>0$
in their maximal domain of definition.  Outputs take values in $U$ as well,
that is, $h:X\rightarrow U$.  

Convex proper cones with nonempty interiors are given on $\R^m$ and $\R^n$,
and, the system is monotone with respect to the orders that they induce.  That
is, $x(0)\leq x'(0)$ and $u(t)\leq u'(t)$ for all $t$ imply $x(t)\leq x'(t)$ for all $t$
in the common domain of definition, if $\dot x(t)=f(x(t),u(t))$ and
$\dot x'(t)=f(x'(t),u'(t))$.

A standard case is that in which $X$ and $U$ are the respective positive
orthants (that is, states and inputs are required to have nonnegative
coordinates), and the cones are also orthants, corresponding to coordinatewise
partial orders, see~\cite{angeli-TAC}.  See~\cite{smith,smith-hirsch} for
excellent introductions to monotone systems (but with no inputs or outputs).

When $h$ is a monotone mapping, that is to say, $x\leq x'\Rightarrow h(x)\leq h(x')$,
we say that the system~(\ref{sys1}) is {\em input/output (i/o) monotone.\/}
If $h$ is anti-monotone, meaning that $x\leq x'\Rightarrow h(x)\geq h(x')$, we will
say that the system~(\ref{sys1}) is {\em input/output (i/o) anti-monotone.\/}
Of course, these definitions are all relative to the orders in states and
outputs. 

{\em We will suppose from now on that the system~(\ref{sys1}) is i/o anti-monotone.\/}

Now we consider the system~(\ref{sys2}).  Its state space is the Cartesian
product $X\times X$.  As an order in this extended state space, we will pick, for
pairs $(x,z)$ and $(x',z')$, the following one:
\[
(x,z) \leq (x',z') \;\;\Leftrightarrow\;\;
x\leq x' \mbox{\ and\ } z\geq z'
\]
where ``$\leq $'' denotes the order on $X$.
The extended system is monotone.   Indeed, suppose that
$(x(0),z(0)) \leq  (x'(0),z'(0))$ in the order in $X\times X$, that is,
$x(0)\leq x'(0)$ and $z'(0)\leq z(0)$, and that $u(t)\leq u'(t)$.
Then $x(t)\leq x'(t)$ for all $t$ (where defined).
Now consider $h(x(t))$ and $h(x'(t))$ as inputs to the system $\dot z=f(z,x)$,
and notice that $h(x'(t))\leq h(x(t))$ since $h$ is anti-monotone.
Since $z'(0)\leq z(0)$, monotonicity of~(\ref{sys1}) gives us
that $z'(t)\leq z(t)$ for all $t$.
Thus, $(x(t),z(t)) \leq  (x'(t),z'(t))$ in the order in $X\times X$,
showing monotonicity.
It is very important to notice that the output map $(x,z)\mapsto h(z)$
is monotone, not anti-monotone, once that we have reversed the order in this
fashion.  
Therefore, {\em the extended system~(\ref{sys2}) is i/o monotone.\/}

We will make the following assumption:
\bi
\item[(B)]
Every solution of the closed loop system $\dot x=f(x,h(z))$, $\dot z=f(z,h(x))$
is bounded.
\ei
Although this is, generally speaking, a strong assumption,
it is automatically satisfied in many situations of interest,
as we discuss below.  Moreover, our purpose is to re-interpret small-gain
theorems as in~\cite{angeli-TAC,german-DCDS} in terms of the extended system,
and small-gain theorems have similar boundedness hypotheses.

One general situation in which (B) is satisfied is as follows.  Suppose that
the open-loop system $\dot x=f(x,u)$ has the following property:
{\em for every constant input $u(t)\equiv u_0$, solutions are bounded.}
One case when this happens is if, as in~\cite{angeli-TAC,german-DCDS}, one
assumes that a {\em characteristic\/} exists, meaning that for every constant
input, there is a state $x_u\in X$ such that all solutions converge to $x_u$.  
Existence of characteristics obviously implies boundedness.
(As it turns out, boundedness together with uniqueness of steady states for
each $u_0$ is already equivalent to the convergence property, via the results
of~\cite{JiFa,Dancer}.)
Suppose also that, as common in biological inhibitory loops, the map
$h$ is bounded (for instance, it has the form $V/(K+x_n)$).
Finally, assume that $X$ and $U$ have the property that any bounded subset
is included in a ``rectangle'' in the sense that $U_0\subseteq U$ bounded implies the
existence of $u_-,u_+\in U$ such that $u_-\leq u\leq u_+$ for all $u\in U_0$ and
similarly for $X$.
Now consider any solution of the closed-loop system.  Since $u=h(z)$
is bounded, this solution is also a solution of the open-loop system
$\dot x=f(x,u)$, $\dot z=f(z,h(x))$ with some bounded input $u$.  So there are two
constant inputs $u_-$ and $u_+$ such that $u_-\leq u(t)\leq u_+$ for all $t$.
This means, because of monotonicity, that $x(t)$ remains bounded between the
solutions corresponding to the two constant inputs.  Now we repeat with the
system $\dot z=f(z,h(x))$, which has a bounded input since $h(x(t))$ is bounded.
So the state of the composite system remains bounded, as we wanted.

\section*{An Observation}

The main observation is as follows:
\begin{quotation}
\noindent
{\em Suppose that the closed loop system $\dot x=f(x,h(z))$, $\dot z=f(z,h(x))$
has a unique steady state.
Then there is a state $\barx$ such that every solution of $\dot x=f(x,h(x))$
converges to $\barx$.}
\end{quotation}
The assumption means that there is only one solution $(\barx,\barz)$ to these
algebraic equations: 
\beqn
f(x,h(z))&=&0\\
f(z,h(x))&=&0 \,.
\eeqn
From the results in~\cite{JiFa,Dancer} (Theorem 1 in the latter,
although stated for discrete iterations, also applies to differential
equations, as remarked in that same paper), we have then that all solutions
of the closed-loop system converge to $(\barx,\barz)$.  (Note that monotonicity
of the extended system is essential in this argument.  The results
from~\cite{JiFa,Dancer} do not apply to the original system, which is not
monotone in closed-loop.)
In particular, then, we can consider those solutions starting at the diagonal
$\Delta =\{(x,x)\in X\times X\}$.  By uniqueness of solutions, these are precisely the
functions of the form $(x(t),x(t))$, where $x(t)$ solves $\dot x=f(x,h(x))$.
Therefore, we conclude that $x(t)\rightarrow \barx$.

We have given a sufficient condition ensuring convergence of all solutions to
an equilibrium,
for the {\em negative feedback\/} (and hence not monotone) system
$\dot x=f(x,h(x))$.
This condition is far from being necessary.
For example, consider the linear monotone system $\dot x=-x+u$, with $X=U=\R$
and the output map $h(x)=-kx$, where $k$ is a positive constant.
The closed-loop system $\dot x=f(x,h(x))$ is $\dot x=-(1+k)x$, which is globally
asymptotically stable for any $k>0$.
On the other hand, the composite system is $\dot \xi =A\xi $, where $\xi $ is the
vector with components $x$ and $z$ and
\[
A = \pmatrix{-1&-k\cr
             -k&-1}
\]
and this matrix is stable (so (B) holds, and equilibria are unique)
if and only if $k<1$ (``small-gain'' condition).
The gap is due to the fact that arbitrary delays can be tolerated in the
feedback loop: the same conditions guarantee that the delay-differential
system $\dot x(t)=f(x(t),h(x(t-r)))$ has a global stability property, for any
delay length $r>0$, as follows by applying the same argument to the extended
system.  (See~\cite{smith} for a discussion of how delays, if appearing at
appropriate places in the system equations, do not affect stability for
monotone systems, or, specifically in the context of small gain theorems the
discussion in~\cite{lsu-book}.)  

We next discuss how our observation relates to the conditions imposed in
small-gain theorems for monotone systems.
For this, we assume that the system $\dot x=f(x,u)$ admits a characteristic,
and we denote as $k:U\rightarrow U$ the function that assigns to each value $u\in U$
the output value $y=h(x_u)$, where $x_u$ is the steady state to
which trajectories converge under the constant input $u(t)\equiv u$.
It is easy to see that $k$ is anti-monotone.  Thus, if there is a fixed point
$k(u)=u$, then under slight additional assumptions (strict decrease), this
fixed point must be unique.  On the other hand, convexity of $U$ plus
appropriate boundedness often allows one to apply a fixed-point theorem and
conclude that such an $u$ must indeed exist.
In any event, the small-gain conditions in~\cite{angeli-TAC,german-DCDS}
amount to asking that the following iteration on $U$:
\[
u^+ = k(u)
\]
must globally converge to a fixed point.
As remarked in this context in~\cite{patrick-MBE} for scalar systems,
and generally in~\cite{german-DCDS},
this global convergence condition is equivalent, under assumptions subsumed by
the ones made here, to asking that solutions of the iteration be bounded (which
automatically happens if, for example, $h$ is bounded) and that
\[
k(k(u))=u \mbox{\ \ \ has a unique solution.}
\]
Now, the characteristic of the extended system~(\ref{sys2}) is well-defined,
and it is precisely $k^2$.
In view of the one to one correspondence between steady states of the
closed-loop extended system and solutions of $k(k(u))=u$, we conclude that
the assumption that the closed loop system $\dot x=f(x,h(z))$, $\dot z=f(z,h(x))$
system has a unique steady state is equivalent to the small-gain condition
that $u^+ = k(u)$ is globally convergent.
This provides an equivalence between the condition as stated here
and the small-gain condition as stated in our previous papers.

\subsection*{Linear Systems}

The picture is especially clear for linear systems.  Linear systems are
of interest in themselves and also serve to provide local versions of
stability results for nonlinear systems.  For linear systems, one can give
a self-contained algebraic proof, with no need to appeal to~\cite{JiFa,Dancer}. 

We assume now that $f(x,u)=Ax+Bu$, $h(x)=-Cx$, where the matrix $A$ is
an $n\times n$ quasi-monotone real matrix and $B$ and $C$ are both monotone,
$B\in \R^{n\times m}$ and $C\in \R^{m\times n}$.  (Quasi-monotonicity means that
the exponentials $e^{tA}$ are monotone, for all $t>0$, and can be
characterized algebraically, see e.g.~\cite{smith-hirsch}.
For the standard case in which the cone is the main orthant in $\R^n$,
quasi-monotone is the same as ``Metzler matrix,'' i.e. the off-diagonal
entries of $A$ must be nonnegative.)
We also assume that $A$ is a Hurwitz matrix (all its eigenvalues have negative
real parts).
The extended open-loop system has the form
$\dot \xi =F\xi +Gu$, $y=H\xi $, where
\[
F = \pmatrix{A &  0\cr
           -BC &  A},
\quad
G = \pmatrix{B\cr
             0},
\quad
H = \pmatrix{0 & -C},
\quad
\]
and its characteristic is the linear map $k^2(x)=K^2x$, where
\[
K = -CA^{-1}B \,.
\]
The closed-loop extended system is $\dot \xi =(F+GH)\xi $, where
\[
F+GH = \pmatrix{A & -BC\cr
              -BC &  A}\,.
\]
It is shown in~\cite{angeli-SCL} for the scalar input case, and
in~\cite{german-SCL} (Theorem 2 and Lemma 1; the nonsingularity assumption
turns out to be redundant) for the general case,
that for a linear monotone system, stability of the closed loop amounts to
the spectral radius $\rho $ of the characteristic being less than one.
That is, $F+GH$ is a Hurwitz matrix if and only if all eigenvalues of
$K^2$ lie strictly inside the unit disk. 

A matrix $K$ has spectral radius less than one if and only if $K^2$ does.
Therefore, $F+GH$ being Hurwitz is equivalent to the small-gain condition
$\rho (K)<1$.  (Which is equivalent, of course, to the global convergence of the
iteration $u^+=k(u)$.)

A change of variables using $y=x+z$ on the extended system brings it into the
following form:
\beqn
\dot y &=& (A-BC) y\\
\dot z &=&    -BC y + (A+BC) z \,,
\eeqn
so stability of the composite system is equivalent to:
\[
A-BC \mbox{\ and\ }  A+BC  \mbox{\ are both Hurwitz matrices.}
\]
The gap with necessity is then clear: stability of the original system under
feedback means merely that $A-BC$ must be Hurwitz.

It is also interesting to notice that, when Theorem 2 and Lemma 1
of~\cite{german-SCL} are applied to the i/o monotone system $\dot x=Ax+Bu$ with
output $y=Cx$ (instead of the anti-monotone output $y=-Cx$), one concludes
that $\rho (K)<1$ is equivalent to just $A+BC$ being a Hurwitz matrix.

\section*{Concluding Remarks}

We described how ``small-gain'' results for negative feedback loops around
monotone systems can be viewed as a consequence of results for feedback loops
on an extended monotone system.  It is important to point out, however, that
often small-gain results can be derived under somewhat weaker assumptions that
those made here, and our ``proof by embedding'' may not cover all such
situations.

Negative feedback loops involving monotone systems are of interest, in
particular, because any system with sign-definite interactions can be written
in that form.  This is obvious (just pull-out the negative connections into
the feedback loop), but can be analyzed from the point of view of minimizing
the number of inputs and outputs (cf.~\cite{german-SCL}).

An interesting consequence of our approach is the following connection between
periodic behavior of the original system and multi-stability of the extended
one:  if the extended system has unique steady states, then we have
convergence to equilibria in the original system.  Thus, periodic behavior of
this system would imply, under the assumptions of this note, the existence of
multiple equilibria for the extension.  

A very weak sort of converse is valid
as well.  Input values $u$ with $k(k(u))=u$ which do not arise from the
unique fixed point of $k(u)=u$ are period-two orbits of the iteration
$u^+=k(u)$.  Now suppose that we consider the delay differential system
$\dot x(t)=f(x(t),h(x(t-r)))$, where the delay $r>0$ is very large.
We take the initial condition $x(t)=x_0$, $t\in [-r,0]$, where $x_0$ is picked
in such a manner that $h(x_0)=u_0$, and $u_0\not= u_1$ are two elements of $U$
such that $k(u_0)=u_1$ and $k(u_1)=u_0$.
By definition of characteristic, the solution $x(t)$ approaches $x_1$,
where $h(x_1)=u_1$, if the input to the open-loop system $\dot x=f(x,u)$
is $u(t)\equiv u_0$.  Thus, if $r$ is large enough, the solution of the closed-loop
system will be close to the constant value $x_1$ for $t\approx r$.
Repeating this procedure, one can show the existence of a lightly damped
``oscillation'' between the values $x_0$ and $x_1$, in the sense of a
trajectory that comes close to these values as many times as desired
(a larger $r$ being in principle required in order to come closer and more
often). 
In certain applications such as molecular biological
experiments, measurements have poor resolution and time duration.  In such a
situation, it may be impossible to practically determine the difference
between such pseudo-oscillations and true oscillations.

It is an open question to prove the existence of true periodic orbits,
for large enough delays, when the small-gain condition fails.  This is closely
related to questions of singular perturbations for delay systems, by
time-re-parametrization (we owe this remark to Roger Nussbaum and to
Hal Smith).  For example, in the scalar case, and with $y=x$,
asking that $\dot x=f(x,x(t-r))$, has periodic orbits for large enough $r$ is
equivalent to asking that $\varepsilon \dot x(t)=f(x(t),x(t-1))$ has periodic orbits for
small enough $\varepsilon >0$.  For $\varepsilon =0$, we have the algebraic equation $f(x,u)=0$
that defines the characteristic $x=k(u)$.  Thus one would want to know that
periodic orbits of the iteration $u^+=k(u)$, seen as the degenerate case
$\varepsilon =0$, survive for small $\varepsilon >0$.

\end{document}